\documentclass[12pt,twoside]{article}
\usepackage{epsfig,amsfonts, amssymb,amsmath}

\newtheorem{Th}{Theorem}[section]
\newtheorem{ThDef}[Th]{Theorem and definition}
\newtheorem{Prop}[Th]{Proposition}
\newtheorem{Lm}[Th]{Lemma}
\newtheorem{Def}[Th]{Definition}

\newtheorem{Rmk}{\sl Remark}[section]

\def\ep{\hfill {$\square$}\\}
\def\bp{{\bf Proof.}\;}
\def\pref#1{(\ref{#1})}
\let\dsp=\displaystyle

\def\r{\mathbb R}
\def\N{\mathbb N}
\def\eps{\varepsilon}
\def\1{{1\hspace{-1.2mm}{\rm I}}}
\def\Lip{\rm Lip}
\def\J{\mathcal J}
\def\W{\mathcal W}
\def\E{\mathcal E}
\def\VV{V^T_e}
\def\V{\mathcal V}
\def\R{\mathcal R}
\def\P{\mathcal P}
\def\S{{\mathcal S}_{BV}}
\def\B{\mathcal B}
\def\diver{\operatorname{div}}
\def\sign{\operatorname{sgn}}
\def\tr{\mathop{\rm tr}}
\def\aa{{\bf a}}
\def\I{\mathop{\rm Id}\nolimits}
\def\diam{\mbox{{\normalsize diam}}}

\def\theequation{\thesection.\arabic{equation}}
\def\thesection{\arabic{section}}
\def\thesubsection{\arabic{section}.\arabic{subsection}}
\def\thesubsubsection{\arabic{section}.\arabic{subsection}.\arabic{subsubsection}}
\newcommand\Section{%
\def\thesubsection{\arabic{section}}
\setcounter{Th}{0}
\setcounter{Rmk}{0}
\setcounter{Ex}{0}
\setcounter{equation}{0}\section}

\begin{document}

\title{Uniqueness and weak stability for multi-dimensional
transport equations with one-sided Lipschitz coefficient}

\author{F. Bouchut$^1$, F. James$^2$, S. Mancini$^3$}
\date{}

\maketitle

{\small\begin{center}
	$^1$ DMA, Ecole Normale Sup\'erieure et CNRS\\
	45 rue d'Ulm\\
	75230 Paris cedex 05, France\\
	e-mail: Francois.Bouchut@ens.fr
\end{center}
\begin{center}
	$^2$ Laboratoire MAPMO, UMR 6628\\
	Universit\'e d'Orl\'eans\\
	45067 Orl\'eans cedex 2, France\\
	e-mail: Francois.James@labomath.univ-orleans.fr
\end{center}
\begin{center}
	$^3$ Laboratoire J.-L. Lions, UMR 7598\\
	Universit\'e Pierre et Marie Curie, BP 187\\
	4 place Jussieu, 75252 Paris cedex 05, France\\
	e-mail: smancini@ann.jussieu.fr
\end{center}}
\vspace{2cm}

\begin{abstract} The Cauchy problem for a multidimensional linear
transport equation with discontinuous coefficient is investigated.
Provided the coefficient satisfies a one-sided Lipschitz condition,
existence, uniqueness and weak stability of solutions are obtained
for either the conservative backward problem or the advective
forward problem by duality.
Specific uniqueness criteria are introduced for the backward
conservation equation since weak solutions are not unique.
A main point is the introduction of a generalized flow in the sense of 
partial differential equations, which is proved to have unique 
jacobian determinant, even though it is itself nonunique.
\end{abstract}



\newpage
\baselineskip=12pt
\parskip=0pt plus 1pt

\tableofcontents
\bigbreak

\Section{Introduction}\label{Intro}
We consider the transport equation
\begin{equation}\label{transport}
	\partial_t u + \aa\cdot \nabla u =0
	\qquad \mbox{in}\ (0,T)\times\r^N,
\end{equation}
with initial data
\begin{equation}
	u(0,x)= u^0(x),
	\label{Cauchy}
\end{equation}
where $\aa=(\aa_i(t,x))_{i=1,\dots,N} \in L^\infty( (0,T)\times \r^N )$
can have discontinuities. The transport equation (\ref{transport}) 
naturally arises with discontinuous coefficient $\aa$ in several applications, together with
the conservation equation
\begin{equation}\label{conserv}
	\partial_t \mu + \diver(\aa\mu) =0
	\qquad \mbox{in}\ (0,T)\times\r^N.
\end{equation}
It is well known that both problems are closely related to the notion of
characteristics, or flow. The flow $X(s,t,x)$, $0\le s,t\le T$,
$x\in\r^N$, is classically defined by the ODE
\begin{equation}\label{ClassicFlow}
	\partial_s X = \aa(s,X(s,t,x)), \qquad X(t,t,x)=x.
\end{equation}
Indeed, when $\aa$ is smooth enough, the flow is uniquely determined by (\ref{ClassicFlow}),
and the solutions $u$ to (\ref{transport}) and $\mu$ to \pref{conserv} are given respectively 
by the classical formul\ae\
$u(t,x)=u^0(X(0,t,x))$, $\mu(t,x)=\det(\partial_xX(0,t,x))\mu^0(X(0,t,x))$. 
The whole theory of characteristics fails if $\aa$ is not smooth:
the flow is no longer uniquely defined, and the notion of solution to (\ref{transport})
or \pref{conserv} has to be reinvestigated.

In the specific case of two-dimensional hamiltonian transport
equations, (\ref{transport}) has been solved with continuous
(non differentiable) coefficients by Bouchut and Desvillettes \cite{BD}, and more recently
with $L^p_{loc}$ coefficients by Hauray \cite{Ha}.
The general well-posedness theory and the connection between the ODE
(\ref{ClassicFlow}) and the PDE (\ref{transport}) was investigated by
DiPerna and Lions \cite{DPL} within the framework of renormalized
solutions, under the assumption that $\aa(t,.)$ lies in
$W^{1,1}_{loc}(\r^N)$ and its distributional divergence belongs to
$L^\infty_{loc}$. The renormalized approach was extended by
Bouchut \cite{B} to the Vlasov equation with BV coefficients, and
very recently, Ambrosio \cite{A} gave the full
generalization in the same context, when $\aa(t,.)\in BV_{loc}$ and
$\diver\aa(t,.)$ is merely $L^1_{loc}$. However, the condition on $\aa$
which arises very often in applications, the so-called
one-sided Lipschitz condition (OSLC)
\begin{equation}\label{OSLC}
	\exists \alpha \in L^1_+(0,T)\; {\ \rm s.t.\ }\;
	\langle\aa(t,y)- \aa(t,x),y-x\rangle \leq \alpha(t) |y-x|^2
\end{equation}
for almost every $(t,x,y)\in (0,T)\times \r^N\times \r^N$, 
implies only an upper bound on $\diver\aa$, thus
$\diver\aa$ is not absolutely continuous with respect to the Lebesgue
measure. It turns out that in this situation the behaviour of the
discontinuities is completely different, and the renormalized approach
is not adapted.

For discontinuous $\aa$, Filippov's theory \cite{Filippov}
gives a generalized definition of a solution $X$ to the ODE
\pref{ClassicFlow} for a merely bounded coefficient $\aa$. Uniqueness,
as well as some stability results, are ensured under the OSLC condition \pref{OSLC}.
Notice that, in contrast with
the classical theory of ODEs, uniqueness only holds for the forward problem, i.e. for
$s>t$.
This theory is used to solve differential equations and inclusions, see \cite{F2,Aubin},
and also to define generalized characteristics in the context of nonlinear conservation
laws, as initiated by Dafermos \cite{D}.

The interest of \pref{OSLC} lies in the fact that it allows to prove
rigorous results in the situation where $\diver\aa$ is a negative
measure. This corresponds intuitively to a compressive situation,
which arises in many applications in nonlinear hyperbolic equations.
The conservation equation (\ref{conserv}) appears indeed in the study of several
degenerate hyperbolic conservation laws where solutions are measures in space,
see e.g. Keyfitz and Kranzer \cite{KK}, LeFloch \cite{PGLF} or
Zheng and Majda \cite{ZM}, and the system of pressureless gases, see \cite{BJpg} and the
references quoted therein. It is also naturally involved in the context of linearization
of conservation laws, see \cite{BJdiff} in the one-dimensional setting, \cite{GodOlaRa,GodOlaRa2}
for an application in fluid mechanics, \cite{JS} in the context of an inverse problem.
Several results are to be quoted here, Poupaud and Rascle \cite{PR} use the Filippov flow to study the
multidimensional equation, Popov and Petrovna \cite{Popov,PopovMD} give several examples of flows
and existence results, both in the one and multidimensional cases.

The idea we use here to solve the Cauchy problem
\pref{transport}-\pref{Cauchy} is based on a previous work of Bouchut
and James \cite{BJ}, and makes use of the so-called {\sl duality
solutions}. The main idea consists in solving the dual (or adjoint)
equation to \pref{transport}-\pref{Cauchy}, which turns out to be a conservative
backward problem
\begin{equation}\label{back}
	\partial_t\pi +\diver(\aa\pi)=0\qquad\mbox{in }(0,T)\times\r^N,
\end{equation}
with final data
\begin{equation}
	\pi(T,x)=\pi^T(x).
	\label{finalrever}
\end{equation}
Indeed a formal computation leads to
\begin{equation}
	\partial_t(u\pi) + \diver(\aa u\pi) = 0,
	\label{dualprod}
\end{equation}
and hence to
\begin{equation}\label{def-duality}
	\frac{d}{dt}\int u(t,x)\pi(t,x)\,dx = 0.
\end{equation}
It is quite classical now that existence for
\pref{back}-\pref{finalrever} leads to uniqueness for the corresponding
direct problem, and this has been successfully used in the nonlinear
context, since Oleinik \cite{Ole}, Conway \cite{C}, Hoff \cite{H}.
However, the lack of uniqueness, and therefore of stability, for
weak solutions to \pref{back}-\pref{finalrever} usually forbids the use
of \pref{def-duality} as a convenient definition for solutions
to the direct problem.

The corner stone in \cite{BJ} was the introduction of the {\sl
reversible solutions} to the backward problem
\pref{back}, thus defining a class for which
existence, uniqueness and weak stability hold. The theory of duality
solutions follows then, leading to existence, uniqueness, weak
stability for both the transport equation \pref{transport} and the
conservation one \pref{conserv}, for one space dimension.

In this paper we generalize this approach to multidimension. However,
since we are only able to establish weakly stable uniqueness criterion
for the conservative backward problem \pref{back}-\pref{finalrever},
and not for the nonconservative one, we are therefore only able to deal
with the forward advective, or nonconservative transport equation. The
weak stability of the conservation equation \pref{conserv} is presently out of reach
in multidimension. In particular, generalized flows in the sense of
PDE, which we call {\sl transport flows}, are not unique. However,
since their jacobian determinant is unique, they allow us to define a
notion of reversible solutions to \pref{back}. The various
characterizations we had at hand in \cite{BJ} (positivity,
renormalization) do not hold here, except if $\pi$ can be expressed as
a linear combination of jacobians of locally Lipschitz functions.
However it turns out that this is enough to define duality solutions to
the direct equation \pref{transport}, and obtain convenient stability
results.

The paper is organized as follows. In Section \ref{coefOSLC}, we
collect and prove several useful results on coefficients satisfying
the OSLC condition \pref{OSLC}. Section \ref{BackProblem} is the core
of the paper, it contains the definition of transport flows and
reversible solutions. Next, duality solutions are defined in Section
\ref{forward}, where existence and uniqueness are proved. Stability
results are gathered in Section \ref{stability}. Section
\ref{nonuniqueness} emphasizes the fact that transport flows are not
necessarily unique, and we prove several technical results about
jacobian determinants in the Appendix.

\Section{Properties of coefficients satisfying the OSLC condition}\label{coefOSLC}

This section is devoted to remarkable properties of coefficients
satisfying the OSLC condition that are used in the paper.
Sharper properties can be found in \cite{AA}.
	\begin{Lm} If a sequence $\aa_n$ is bounded in $L^\infty((0,T)\times\r^N)$
and satisfies the OSLC condition with some $\alpha_n(t)$ uniformly bounded
in $L^1(0,T)$, then, up to a subsequence, $\aa_n$ converges in
$L^\infty-w*$ to $\aa\in L^\infty((0,T)\times\r^N)$ satisfying the
OSLC condition with some $\alpha(t)$ such that $\int_0^T\alpha(t)dt
\le\lim\inf \int_0^T\alpha_n(t)dt$.
	\label{Lemma suitean}\end{Lm}
\bp
We can pass to the limit in \pref{OSLC} in the sense of distributions
in $(0,T)\times\r^N\times\r^N$. Then the right-hand side involves some
$\hat\alpha(t)\in {\cal M}_+(0,T)$ such that $\alpha_n\rightharpoonup\hat\alpha$.
We deduce that the required OSLC inequality holds with $\alpha$
the absolutely continuous part of $\hat\alpha$.
\ep

\noindent The OSLC condition implies that $\diver\aa$ is bounded from above,
\begin{equation}
	\diver \aa\leq N\alpha(t).
\end{equation}
This is actually a consequence of the following stronger result.
\begin{Lm}\label{Lemma OSLC}
(i) The OSLC condition \pref{OSLC} is equivalent to
\begin{equation}\label{symm-matr1}
	\frac{\nabla \aa + \nabla \aa^t}{2} \leq \alpha(t)\I
\end{equation}
for {\sl a.e.} $t$, in the sense of matrix distributions in $\r^N$,
which means that
\begin{equation}\label{symm-matr2}
	\forall h\in\r^N\qquad
	\nabla \aa \cdot h \cdot h \leq \alpha (t) |h|^2
\end{equation}
in the sense of distributions in $\r^N$.

\noindent (ii) A coefficient $\aa\in L^\infty$
satisfying the OSLC condition verifies that for a.e. $t$,
$\nabla\aa+\nabla\aa^t$ is locally a matrix valued bounded measure
in $x$. Moreover, for any bounded convex open subset $C$ of $\r^N$,
\begin{equation}
	\forall h\in\r^N\qquad\int_C |\nabla\aa\cdot h\cdot h|
	\le 2\Bigl(\alpha(t)|C|+\|\aa\|_\infty\diam(C)^{N-1}\Bigr)|h|^2.
	\label{bornemes}
\end{equation}
\end{Lm}
\bp We can consider a fixed time $t$.
For (i), assume first that $\aa$ is smooth in $x$.
Then the OSLC condition \pref{OSLC} is equivalent to
\begin{equation}
	\langle\aa(x+\eps h)- \aa(x), \eps h\rangle \leq \alpha |\eps h|^2
	\label{epsh}
\end{equation}
for any $x$, $h\in\r^N$ and $\varepsilon>0$.
Now, since
\begin{equation}
	\langle\aa(x+\eps h)-\aa(x),\eps h\rangle
	=\int_0^1 \langle\nabla \aa(x+\theta\eps h )
	\,\eps h,\eps h\rangle\, d\theta,
	\label{eqgrad}
\end{equation}
if $\nabla \aa(x) \cdot h\cdot h\leq \alpha |h|^2$ for all $x,h$,
we get obviously \pref{epsh}. Conversely, if \pref{epsh} holds, then
\begin{equation}
	\int_0^1 \langle \nabla \aa (x+\theta\eps h)
	h,h\rangle\, d\theta \leq \alpha |h|^2,
	\label{eqintgrad}
\end{equation}
and by letting $\eps\rightarrow 0$ we recover
$\nabla \aa(x) \cdot h\cdot h\leq \alpha |h|^2$.

\noindent For the general case $\aa\in L^\infty$,
consider the convolution of $\aa$ by a smoothing sequence
$\rho_\eps(x)$, $\aa_\eps =\rho_\eps*\aa$, so that $\aa_\eps
\to \aa$ in $L^1_{loc}(\r^N)$.
If $\aa$ satisfies \pref{OSLC}, then since
\begin{equation}
	\langle \aa_\eps(y)-\aa_\eps(x),y-x\rangle
	=\int\langle \aa(y-z)-\aa(x-z),y-z-(x-z)\rangle\rho_\eps(z)\,dz,
	\label{OSLCconvol}
\end{equation}
$\aa_\eps$ also satisfies \pref{OSLC} for all $\eps>0$, with the same
$\alpha$. Thus by the above proof, $\aa_\eps$ satisfies
$\nabla\aa_\eps\cdot h\cdot h\le\alpha|h|^2$, and by letting
$\eps\to 0$, this gives \pref{symm-matr2}.
Conversely, if \pref{symm-matr2} holds, then
$\nabla\aa_\eps\cdot h\cdot h=\rho_\eps*(\nabla\aa\cdot h\cdot h)
\le \alpha |h|^2$ thus by the proof above, $\aa_\eps$ satisfies \pref{OSLC}.
By letting $\eps\to 0$, we finally get that $\aa$ itself
satisfies \pref{OSLC}, and this concludes the proof of (i).

For (ii), we have by \pref{symm-matr2} that $\alpha|h|^2-\nabla\aa\cdot h\cdot h\ge 0$,
thus it is locally a measure. We conclude that $\nabla\aa\cdot h\cdot h$
is locally a measure for any $h$, and by finite linear combinations
that $\nabla\aa+\nabla\aa^t$ is locally a measure.
It remains to prove \pref{bornemes}.
Writing $\nabla\aa\cdot h\cdot h=\alpha|h|^2-\bigl(\alpha|h|^2-\nabla\aa\cdot h\cdot h\bigr)$,
we get
\begin{equation}
	|\nabla\aa\cdot h\cdot h|\le \alpha|h|^2+\bigl(\alpha|h|^2-\nabla\aa\cdot h\cdot h\bigr).
	\label{inegabs}
\end{equation}
The indicator function of $C$ can be approximated by a function
$\varphi(x)$, $\varphi\in C^\infty_c(C)$, $0\le\varphi\le 1$. Then
\begin{equation}\begin{array}{l}
	\dsp\hphantom{= } \int\bigl(\alpha|h|^2-\nabla\aa\cdot h\cdot h\bigr)\varphi\\
	\dsp =\alpha|h|^2\int\varphi+\int \aa(x)\cdot h\,\,\nabla\varphi(x)\cdot h\\
	\dsp \le \alpha |h|^2|C|+\|\aa\|_\infty |h|\int|\nabla\varphi(x)\cdot h|.
	\label{estvar}
	\end{array}
\end{equation}
But the convexity of $C$ enables to find $\varphi$ such that for any
direction $h$, $\int|\nabla\varphi\cdot h|\le 2|h|\diam(C)^{N-1}$,
and this concludes (ii). \ep
\begin{Rmk} \rm By Lemma 3.4 in \cite{B}, \pref{bornemes} is equivalent
to the fact that for any $\omega\subset\subset C$ and any $h\in\r^N$
such that $\overline{\omega}+\overline{B}(0,|h|)\subset C$,
\begin{equation}
	\int_\omega\,\Bigl|\langle\aa(t,x+h)-\aa(t,x), h\rangle\Bigr|\,dx
	\le 2\Bigl(\alpha(t)|C|+\|\aa\|_\infty\diam(C)^{N-1}\Bigr)|h|^2.
	\label{bornetransl}
\end{equation}
\label{Rmk transl}\end{Rmk}
\begin{Rmk} \rm A coefficient $\aa\in L^\infty$ satisfying the OSLC
condition verifies the local inequalities
\begin{equation}\begin{array}{c}
	\dsp 0\le\alpha(t)-\partial_i \aa_i\le N\alpha(t)-\diver \aa,
	\qquad i=1,\dots,N,\\
	\dsp\left|\frac{\partial_j \aa_{i}+\partial_i\aa_{j}}{2}\right|
	\le N\alpha(t)-\diver \aa,\qquad i,j=1,\dots,N,\ i\not=j,
	\label{localest}
	\end{array}
\end{equation}
in the sense of measures. This is because by \pref{symm-matr1},
the symmetric matrix valued measure
$B=\alpha(t)\I-(\nabla\aa+\nabla\aa^t)/2$ is nonnegative,
and thus satisfies $0\le B_{ii}\le \tr B$, $|B_{ij}|\le \tr B$
(for the latter inequality, apply the Cauchy-Schwarz inequality
to a regularization $B^\varepsilon$ of $B$,
$|B^\varepsilon_{ij}|\le\sqrt{B^\varepsilon_{ii}B^\varepsilon_{jj}}$).
Another useful inequality is obtained by writing that $B\le (\tr B)\I$.
Replacing $B$ by its value gives
$\alpha\I-(\nabla\aa + \nabla\aa^t)/2\le (N\alpha-\diver\aa)\I$, or
\begin{equation}
	(\diver\aa)\I-\frac{\nabla\aa + \nabla\aa^t}{2}
	\le(N-1)\,\alpha(t)\I.
	\label{borna}
\end{equation}
Notice in particular that the inequalities \pref{localest} imply that
\begin{equation}\begin{array}{c}
	\diver\aa\in L^1((0,T),L^p_{loc}(\r^N)) \quad\mbox{for some } 1\le p\le\infty \\
\Longrightarrow \qquad
	\nabla\aa+\nabla\aa^t\in L^1((0,T),L^p_{loc}(\r^N)).
\end{array}\end{equation}
In the case $1<p<\infty$, we have even
$\nabla\aa\in L^1((0,T),L^p_{loc}(\r^N))$, by elliptic regularity
since $\Delta \aa_i=\sum_j\partial_j(\partial_j\aa_i+\partial_i\aa_j)
-\partial_i\diver\aa$.
\label{Rmk rega}\end{Rmk}

We end up this section by pointing out in a simplified context the fact
that the OSLC condition \pref{OSLC} is in some sense orthogonal to the
condition of absolute continuity of the divergence of $\aa$
that is involved in the renormalized theory, and in particular
in the work of Ambrosio \cite{A}.

Consider a coefficient $\aa\in L^\infty$ such that
for almost every $t>0$, $\aa(t,.)\in BV_{loc}(\r^N)$, and the total
variation is integrable with respect to time (indeed according
to \cite{AA}, the OSLC condition implies this regularity).
Then, fixing the time,
the matrix-valued measure $\nabla\aa$ can be
decomposed as $\nabla\aa = \nabla^a\aa + \nabla^s\aa$, where
$\nabla^a\aa$ is absolutely continuous and $\nabla^s\aa$ is singular,
with respect to the Lebesgue measure. Denoting by $|\lambda|$ the total
variation of a (matrix-valued) measure $\lambda$, one has the
polar decomposition $\lambda=M|\lambda|$,
with $|M(x)|=1$ for $|\lambda|$-a.e. $x\in\r^N$. We apply this to the
derivative $\nabla\aa$ of $\aa$, and
make use of Alberti's rank one theorem, which asserts that the
corresponding $M$ is rank-one, so that
\begin{equation}
	\nabla\aa = \eta(x)\otimes\xi(x)|\nabla\aa|,
	\label{alberti}
\end{equation}
where $\xi$ and $\eta$ are unit vectors in $\r^N$. The distributional
divergence of $\aa$ is therefore given by $\diver\aa =
\langle\xi,\eta\rangle|\nabla\aa|$, so that $\diver\aa$ is absolutely
continuous with respect to the Lebesgue measure if and only if $\xi$
and $\eta$ are orthogonal $|\nabla^s\aa|$-a.e. As it is recalled by
Ambrosio in \cite{A}, these results hold true in a very general
context, but we wish to consider merely the case where the jump set of
the function $\aa$ is (locally) a hypersurface in $\r^N$. Then the
vector $\xi$ turns out to be the normal to the hypersurface, while
$\eta$ is the direction of the jump of $\aa$. Thus the usual assumption
of absolute continuity of the divergence implies that the jump of $\aa$
is tangent to the hypersurface.

If now $\aa$ satisfies the OSLC condition (instead of having
absolutely continuous divergence), then the conclusion is radically
different, as the following computation shows.
Consider a point $x$ in the hypersurface and a neighborhood $\omega\ni
x$. The hypersurface divides $\omega$ in two subsets $\omega_-$ and
$\omega_+$, and we take the normal $\xi(x)$ pointing to $\omega_+$. We
assume that the limits $\aa_+(x)=\lim_{y\in\omega_+}\aa(y)$ and
$\aa_-(x)=\lim_{y\in\omega_-}\aa(y)$ exist. We rewrite the OSLC
condition \pref{OSLC} by letting $x$ go to the hypersurface in
$\omega_-$,
and by writing $y=x+\epsilon h$, for any unit vector $h$ and $\epsilon>0$.
After division by $\epsilon$, we get
\begin{equation}
	\langle \aa(x+\epsilon h)-\aa_-(x),h \rangle \le \alpha\epsilon.
	\label{oslcxy}
\end{equation}
Now, if $\langle\xi,h\rangle > 0$, $x+\epsilon h\in\omega_+$ for
$\epsilon$ small enough, thus letting $\epsilon$ go to zero, we obtain
\begin{equation}
	\langle \aa_+(x)-\aa_-(x),h \rangle \equiv \langle [\aa],h\rangle \le 0
	\quad\mbox{for any $h$ such that }|h|=1,\ \langle\xi,h\rangle>0.
	\label{proplimoslc}
\end{equation}
This inequality holds true for any unit vector $h$ such that
$\langle\xi,h\rangle = 0$ as well, so that $[\aa]$ has to be colinear
with the normal $\xi$. Finally, using \pref{proplimoslc} again
we obtain
\begin{equation}
	[\aa]=-\lambda\xi\mbox{ for some }\lambda\ge 0.
	\label{jumpaxi}
\end{equation}
This is radically different from the case above of absolute continuity
of the divergence, where $[\aa]\cdot\xi=0$.

\Section{Backward problem, reversible solutions}\label{BackProblem}

In this section we intend to study weak solutions
$\pi\in C([0,T],L^\infty_{loc}(\r^N)w*)$ to the conservative
transport equation \pref{back}, together with solutions $p\in \Lip_{loc}([0,T]\times\r^N)$ to
the nonconservative transport equation
	\begin{equation}\label{noncons}
\partial_t p + \aa\cdot \nabla p =0\qquad \mbox{a.e. in }(0,T)\times\r^N.
	\end{equation}
We prove very general properties of these solutions, and indeed
all that is done in this section is true for a coefficient
$\aa\in L^\infty$ such that there exists a Lipschitz transport flow
(see Definition \ref{Def flow}). Only the existence result
Proposition \ref{Prop existflow} makes use of the OSLC condition.\\

We denote by $\W^T$ the vector space of all
$\pi\in C([0,T],L^\infty_{loc}(\r^N)w*)$ solving \pref{back}
in the weak sense. Recall that there is no uniqueness for the Cauchy
problem in $\W^T$, so that we introduce various restrictive definitions
of convenient solutions. Thus we define several subspaces of $\W^T$,
and the main goal of this section is to establish inclusions between
these subspaces. As we shall see, the definitions do not coincide in
the general case.

The first two notions of solution, and the most natural ones, are defined
through sign or cancellation properties.
We need to introduce the space $\E^T$ of exceptional solutions
to the nonconservative equation,
\begin{equation}
	\E^T = \{ p \in \Lip_{loc}([0,T]\times\r^N)
	\mbox{ solving \pref{noncons} such that} \; p(T,.)=0 \},
	\label{E}
\end{equation}
and the corresponding open support
\begin{equation}
	\VV = \bigcup_{p_e\in \E^T} \{ (t,x)\in (0,T)\times\r^N
	\ \mbox{such that}\ p_e(t,x) \neq 0 \}.
	\label{V}
\end{equation}
Our first notion of ``good'' solution to the conservative equation
\pref{back} is defined via the cancellation in $\VV$,
\begin{equation}
	\V^T = \{ \pi \in \W^T \ \mbox{such that} \
	\pi =0 \ {\rm a.e.\ in}\ \VV \}.
	\label{VT}
\end{equation}
Nonnegative solutions will also play a natural part here, thus we define
\begin{equation}
	\P^T=\{\pi\in\W^T\,;\, \pi\ge 0\}.
	\label{PT}
\end{equation}
Before giving other definitions, we prove the following preliminary result.
	\begin{Prop}\label{Lemma posVT}
We have $\mbox{Vect}(\P^T) \subset \V^T$, that is any linear combination of nonnegative solutions
to \pref{back} satisfies the cancellation property \pref{VT}.
	\end{Prop}
The proof relies on the following lemma, which defines some kind of weak product.
	\begin{Lm}\label{Lemma weakprod}
Let $\pi\in\W^T$ and $p \in \Lip_{loc}([0,T]\times\r^N)$
solve \pref{noncons}. Then $p\,\pi\in\W^T$.
	\end{Lm}
\bp
Let $\pi_\varepsilon$ and $(\aa\pi)_\varepsilon$ denote the
approximations of $\pi$ and $\aa\pi$ by convolution by a smoothing
kernel in $(t,x)$. Since
$\partial_t\pi_\varepsilon+\diver(\aa\pi)_\varepsilon=0$
and using the equation satisfied by $p$ we have
\begin{equation}\begin{array}{l}
	\dsp \partial_t(p\pi_\varepsilon)+\diver\left(p(\aa\pi)_\varepsilon\right)
	=\pi_\varepsilon\partial_t p+(\aa\pi)_\varepsilon\cdot\nabla p\\
	\dsp\hphantom{\partial_t(p\pi_\varepsilon)+\diver\left(p(\aa\pi)_\varepsilon\right)}
	=\left((\aa\pi)_\varepsilon-\aa\pi_\varepsilon\right)\cdot\nabla p,
	\label{calpnpi}
	\end{array}
\end{equation}
thus letting $\varepsilon\rightarrow 0$ we get
\begin{equation}
	\partial_t(p\pi)+\diver\left(\aa \,p\pi\right)=0
	\qquad\mbox{in }(0,T)\times\r^N.
	\label{eqlim}
\end{equation}
Since $p\pi\in C([0,T],L^\infty_{loc}(\r^N)w*)$, we have the result.\ep

\noindent{\bf Proof of Proposition \ref{Lemma posVT}}\;
Let $\pi\in\W^T$ such that $\pi\ge 0$. We actually prove that for any
any $p_e\in\E^T$,
we have $p_e\,\pi=0$ almost everywhere in $(0,T)\times\r^N$.
This obviously implies the result.
Since $|p_e|$ also belongs to $\E^T$ when $p_e\in\E^T$,
we can assume without loss of generality that $p_e\ge 0$.
Since $p_e\pi\in C([0,T],L^\infty_{loc}(\r^N)w*)$, $p_e\pi\ge 0$,
and $(p_e\pi)(T,.)=0$, we conclude with Lemma \ref{Lemma weakprod} and
by integration on cones that $p_e\pi\equiv 0$.
The conclusion that $\pi\in\V^T$ follows obviously by the definition
of $\V^T$, since each set in the right-hand side of \pref{V} is open.\ep

We give now two new definitions of weak solutions to \pref{back}. The first one takes advantage of
an additional conservation law satisfied by the jacobian matrix of a set of $N$ solutions to the
transport problem \pref{noncons}. Actually such a property is quite classical in the context of
elastodynamics: the jacobian and cofactor matrices of a velocity field satisfy
conservation equations, see for instance Quin \cite{Quin}, or Wagner \cite{wagner} for a survey in
the context of conservation laws.
	\begin{ThDef}[Jacobian solutions]
Let $p_i\in \Lip_{loc}([0,T]\times\r^N)$, $1\le i\le N$, solve
	\begin{equation}\label{system_p}
\partial_t p_i + \aa\cdot \nabla p_i =0\qquad \mbox{a.e. in }(0,T)\times\r^N.
	\end{equation}
Then the jacobian determinant
	\begin{equation}
\pi=J(p_1,\dots,p_N)\equiv\det(\nabla p_1,\ldots,\nabla p_N)
\in C([0,T],L^\infty_{loc}(\r^N)w*)
	\label{DefJaco}\end{equation}
is a weak solution to \pref{back}.
We shall call such a solution $\pi$ to \pref{back} a {\sl jacobian
solution}, and denote by $\J^T\subset\W^T$ the set of these solutions.
	\label{Lemma jacobian}\end{ThDef}
\bp
Let us apply Lemma \ref{Lemma divnul} in the Appendix, in the variables
$(t,x)\in (0,T)\times\r^N\subset\r^{N+1}$, to the hamiltonians
$H_i=p_i$, $i=1,\dots,N$. We get that $\diver_{t,x} V=0$,
where $V=(V_0,V_1,\dots,V_N)$ has to be computed by \pref{Ik}.
We have $V_0=\det(\partial p_i/\partial x_j)$, $1\le i,j\le N$,
i.e. $V_0=\pi$. Then, for $1\le k\le N$, denoting $p=(p_1,\dots,p_N)$,
\begin{equation}\begin{array}{l}
	\dsp V_k=(-1)^k\det\left(\partial_tp,\frac{\partial p}{\partial x_1},
	\dots,\widehat{\frac{\partial p}{\partial x_{k}}},
	\dots,\frac{\partial p}{\partial x_{N}}\right)\\
	\dsp\hphantom{V_k}=(-1)^k\det\left(-\sum_{j=1}^N\aa_j\frac{\partial p}{\partial x_j},\frac{\partial p}{\partial x_1},
	\dots,\widehat{\frac{\partial p}{\partial x_{k}}},
	\dots,\frac{\partial p}{\partial x_{N}}\right)\\
	\dsp\hphantom{V_k}=-\aa_k(-1)^k\det\left(\frac{\partial p}{\partial x_k},\frac{\partial p}{\partial x_1},
	\dots,\widehat{\frac{\partial p}{\partial x_{k}}},
	\dots,\frac{\partial p}{\partial x_{N}}\right)\\
	\dsp\hphantom{V_k}=\aa_k\,\pi,
	\label{Vkk}
	\end{array}
\end{equation}
which gives \pref{back}. The time continuity of $\pi$ comes
from the weak continuity of the jacobian, see Lemma \ref{Lemma weakstab}.
Indeed, if $t_n\rightarrow t$,
then $p(t_n,.)\rightarrow p(t,.)$ locally uniformly with uniform
bound in $\Lip_{loc}$ thus $J(p(t_n,.))\rightharpoonup J(p(t,.))$.\ep

We have some kind of uniqueness for Jacobian solutions, in the following sense.
	\begin{Lm}\label{Lemma uniqJ}
Let $p_i,q_i\in \Lip_{loc}([0,T]\times\r^N)$, for $i=1,\ldots, N$,
solve \pref{system_p}, with $p_i(T,.) = q_i(T,.)$.
Recall that this implies by Theorem \ref{Lemma jacobian} that
$J(p)$, $J(q)\in\W^T$.
Assume moreover that $J(p)$, $J(q)\in\V^T$
(by Proposition \ref{Lemma posVT} this is the case in particular
if $J(p)\ge 0$ and $J(q)\ge 0$). Then $J(p)=J(q)$ almost everywhere.
	\end{Lm}
\bp On the one hand, $u_i\equiv p_i-q_i$ is a
locally Lipschitz solution to the transport equation \pref{noncons},
and $u_i(T,.)=0$, i.e. $u_i\in\E^T$.
On the other hand, if we define $v=J(p)$ or $v=J(q)$, we have
$v\in\V^T$, and by definition, we get $v=0$ a.e. in $\{u_i\not =0\}$.
Now, on the set where $p\ne q$,
this gives $v=0$ a.e., thus $J(p)=J(q)=0$. On the set where $p=q$, we
have $\nabla p=\nabla q$ a.e. (see for example \cite{B}, Theorem 2.1),
so that $J(p)=J(q)$ as well.\ep

The second and last notion of solution we give is the most convenient one from the viewpoint
of uniqueness and stability results. However, it does not enjoy really handable characterizations.
It is defined through a specific notion of flow, which we introduce now.
	\begin{Def}[Transport flow]
We say that a vector $X^T\in\Lip([0,T]\times\r^N)$
is a transport flow if
\begin{equation}
	\partial_t X^T+\aa\cdot\nabla X^T=0
	\quad\mbox{in }(0,T)\times\r^N,
	\qquad X^T(T,x)=x,
	\label{flow}
\end{equation}
and $J(X^T)\ge 0$.
	\label{Def flow}\end{Def}

	\begin{Rmk} \rm The definition requires that $X^T$ is globally
Lipschitz continuous. This condition ensures that
$|X^T(t,x)-x|\le \|\partial_t X^T\|_{L^\infty}(T-t)$, and therefore
that $X^T(t,.)$ tends to infinity at infinity.
This estimate can be interpreted as a kind of finite speed
of propagation property.
	\label{Rmk proper}\end{Rmk}

\noindent We first provide a sufficient condition for the existence of
a transport flow.
	\begin{Prop} For any coefficient $\aa\in L^\infty( (0,T)\times \r^N )$
satisfying the OSLC condition, there exists a transport flow $X^T$.
	\label{Prop existflow}\end{Prop}
\bp
Consider a sequence of coefficients $\aa_n\in C^1$, $\aa_n$
bounded in $L^\infty$, satisfying the OSLC condition for some $\alpha_n$
bounded in $L^1$, such that $\aa_n\to \aa$ in $L^1_{loc}$.
Then we have a transport flow associated to $\aa_n$,
$X_n^T(t,x)=X_n(T,t,x)$ where $X_n$ is the classical flow,
\begin{equation}
	\partial_sX_n(s,t,x)=\aa_n(s,X_n(s,t,x)),
	\qquad X_n(t,t,x)=x.
	\label{Xn}
\end{equation}
According to the OSLC condition, we have
\begin{equation}\begin{array}{l}
	\dsp\hphantom{=\ } \partial_s| X_n(s,t,y)-X_n(s,t,x)|^2/2\\
	\dsp=\left\langle\aa_n(s,X_n(s,t,y))-\aa_n(s,X_n(s,t,x)),
	X_n(s,t,y)-X_n(s,t,x)\right\rangle\\
	\dsp\le \alpha_n(s)|X_n(s,t,y)-X_n(s,t,x)|^2,
	\label{accXn}
	\end{array}
\end{equation}
thus for $0\le t\le s\le T$,
\begin{equation}
	| X_n(s,t,y)-X_n(s,t,x)|\le |y-x|
	e^{\int_t^s\alpha_n(\tau)d\tau},
	\label{accXnbis}
\end{equation}
which yields a uniform Lipschitz estimate in $x$ for $X_n$,
\begin{equation}
	\Lip(X_n(s,t,.))\le e^{\int_t^s\alpha_n(\tau)d\tau}.
	\label{lipXn}
\end{equation}
But \pref{Xn} yields a Lipschitz estimate in $s$, while the
transport equation $\partial_tX_n+\aa_n\cdot\nabla X_n=0$
yields a Lipschitz estimate in $t$.
Therefore, up to a subsequence, $X_{n} \to X$ locally uniformly.
We have that $X(s,t,x)$ is Lipschitz continuous in the domain
$0\le t\le s\le T$, $x\in\r^N$, $\partial_tX+\aa\cdot\nabla X=0$
and $X(t,t,x)=x$. By weak continuity of the jacobian
$J(X_n) \rightharpoonup J(X)$ in $L^\infty -w*$, thus $J(X)\ge 0$.
We conclude that $X^T(t,x)=X(T,t,x)$ is a transport flow.
\ep

We emphasize that, in general, there is no uniqueness in the notion of transport flow,
an explicit counterexample is provided in Section \ref{nonuniqueness} below.
However, we have the following uniqueness result for the Jacobian.
	\begin{Prop}
All possible transport flows $X^T$ have the same Jacobian determinant $J(X^T)$.
	\label{Unik jaco flow}\end{Prop}
	\bp If $X^T$ is a transport flow, then $J(X^T)\geq 0$, thus
by Theorem \ref{Lemma jacobian} and Proposition \ref{Lemma posVT}, $J(X^T)\in\V^T$. Therefore Lemma \ref{Lemma uniqJ} ensures
that $J(X^T)$ is uniquely determined. More precisely, $X^T$ is uniquely defined outside $\VV$,
thus $J(X^T)$ also. Now, according to Proposition \ref{Lemma posVT},
we have $J(X^T)=0$ in $\VV$, so that indeed $J(X^T)$ is uniquely determined everywhere. \ep

We are now in position to define the final notion of "good"
solutions to \pref{back}.
	\begin{Def}[Reversible solutions]\label{rev-sol}
We say that $\pi\in\W^T$ is a {\it reversible} solution to \pref{back},
and we shall denote $\pi\in\R^T$,
if for some transport flow $X^T$ one has
\begin{equation}
	\pi(t,x)=\pi(T,X^T(t,x))J(X^T)(t,x)
	\label{RT}
\end{equation}
for all $0\le t\le T$ and a.e. $x\in\r^N$.
	\end{Def}
We notice that this formula is meaningful.
Indeed, according to Lemma \ref{Lemma detnulae}, at a fixed time $t$,
and for any Borel set $Z\subset \r^N$ such that $|Z|=0$,
we have that $J(X^T)=0$ a.e. where $X^T\in Z$.
This proves, together with Remark \ref{Rmk proper}, that the right-hand side of \pref{RT} is well-defined
as an element of $L^\infty_{loc}(\r^N)$ for each $t$.
Also, according to Proposition \ref{Unik jaco flow}, if \pref{RT} holds
for some transport flow $X^T$, then it holds for all transport flows.
Therefore, $\R^T$ is a vector space.
\begin{Th}[Conservative backward Cauchy problem]\label{sol_rev}
For any $\pi^T\in L^\infty_{loc}(\r^N)$,
there exists a unique reversible solution $\pi\in\R^T$ to \pref{back}
such that $\pi(T,.)=\pi^T$.
\end{Th}
\bp
Uniqueness is obvious from the definition.
We claim that a reversible solution with final data
$\pi^T\in L^\infty_{loc}(\r^N)$ can be obtained by the formula
\begin{equation}
	\pi=\pi^T(X^T)J(X^T),
	\label{piX}
\end{equation}
where $X^T$ is any transport flow.
Indeed, for the same reason as above, $\pi(t,.)$ is well-defined
as an element of $L^\infty_{loc}(\r^N)$ for each $t$.
This definition does not depend on the choice of $X^T$, and
obviously $\pi(T,.)=\pi^T$. It only remains to prove that $\pi$
is time continuous and satisfies \pref{back} in the weak sense.

Assume first that $\pi^T$ is locally Lipschitz continuous. Then
$p\equiv\pi^T(X^T)\in\Lip_{loc}([0,T]\times\r^N)$ solves
$\partial_t p+\aa\cdot\nabla p=0$, and by Lemma \ref{Lemma weakprod}
we get that $\pi=p\,J(X^T)\in\W^T$.

In the general case $\pi^T\in L^\infty_{loc}$, we can approximate
$\pi^T$ by $\pi^T_n\in C^1$ such that $\pi^T_n$ is bounded in
$L^\infty_{loc}$ and $\pi^T_n\rightarrow \pi^T$ a.e. Let $Z\subset\r^N$
be the set where $\pi^T_n(x)$ does not converge to $\pi^T(x)$,
with $|Z|=0$. Then for a fixed $t$, $\pi_n(t,x)\equiv\pi^T_n(X^T(t,x))J(X^T)(t,x)$
$\rightarrow \pi^T(X^T(t,x))J(X^T)(t,x)=\pi(t,x)$ a.e., because
$J(X^T)(t,x)=0$ a.e. where $X^T(t,x)\in Z$.
Then, the equation
$\partial_t\pi_n+\diver(\aa \pi_n)=0$ ensures uniform time continuity,
thus $\pi_n\rightarrow\pi$ in $C([0,T],L^\infty_{loc}(\r^N)w*)$
and $\pi$ is a weak solution to \pref{back}.\ep

We look now for intrinsic conditions characterizing
the reversible solutions, that do not involve any transport flow.
It follows from formula \pref{piX} and from Proposition \ref{Lemma posVT} that
\begin{equation}
	\R^T\subset \mbox{Vect}(\P^T)\subset\V^T.
	\label{inclusion}
\end{equation}
In the one space dimension case, these inclusions are equalities, see
\cite{BJ}. Reversibility can also be characterized in terms of
renormalization, namely $\pi$ is reversible if and only of $|\pi|$ is
a weak solution, and total variation properties.
In the present multidimensional context, we are not able to prove the same results,
except for the subclass of linear combinations of jacobian solutions.
This is the aim of the following three propositions.
Note that this is a difficult problem to decompose a bounded
function as a linear combination of jacobian determinants
of lipschitz functions, see \cite{McM}.
	\begin{Prop} Let $\pi\in\W^T$, such that $\pi(T,.)$ can be written
as a linear combination of jacobians of locally Lipschitz functions.
Then $\pi$ is reversible if and only if
\begin{equation}
	\pi\in\mbox{Vect}(\J^T)\cap\V^T.
	\label{jacVT}
\end{equation}
\label{Prop jacVT}\end{Prop}
\bp By assumption, we can write $\pi(T,.)=\sum_{k=1}^K\lambda_kJ(p_k^T)$,
where $p_k^T\in (\Lip_{loc}(\r^N))^N$.
First, if $\pi$ is reversible, then \pref{inclusion} gives
$\pi\in\V^T$, and denoting by $X^T$ a transport flow,
\begin{equation}\begin{array}{l}
	\dsp \pi(t,.)=\pi(T,X^T(t,.))J(X^T(t,.))\\
	\dsp\hphantom{\pi(t,.)}
	=\sum_{k=1}^K\lambda_k\, J(p_k^T)(X^T(t,.))J(X^T(t,.))\\
	\dsp\hphantom{\pi(t,.)}
	=\sum_{k=1}^K\lambda_k\,J\Bigl(p_k^T(X^T(t,.))\Bigr).
	\label{pijac}
	\end{array}
\end{equation}
Since $p_k\equiv p_k^T(X^T)$ is a locally Lipschitz solution to the
nonconservative equation, we deduce that
$\pi\in\mbox{Vect}(\J^T)\cap\V^T$.

\noindent Conversely, let us assume that
$\pi\in\mbox{Vect}(\J^T)\cap\V^T$. We have
\begin{equation}
	\pi=\sum_{k=1}^K\lambda_k\,J(p_k),
	\label{sumjac}
\end{equation}
where $p_k\in(\Lip_{loc}([0,T]\times\r^N))^N$
solves $\partial_tp_k+\aa\cdot\nabla p_k=0$.
Consider the unique reversible
solution $\pi_k^{rev}$ to \pref{back} with final data $J(p_k(T,.))$.
In other words, $\pi_k^{rev}=J(p_k(T,.))(X^T)J(X^T)
=J(q_k)$, with $q_k=p_k(T,X^T)$. Define
$\pi^{rev}=\sum_k\lambda_k\pi_k^{rev}$. Then, since $p_k-q_k\in\E^T$,
for $(t,x)\notin \VV$
we have $p_k(t,x)=q_k(t,x)$ thus $\nabla p_k=\nabla q_k$
a.e. in $(\VV)^c$ and $J(p_k)=J(q_k)$ a.e. in $(\VV)^c$ for all $k$.
We conclude that $\pi(t,x)=\pi^{rev}(t,x)$ a.e. in $(\VV)^c$.
But since $\pi,\pi^{rev}\in\V^T$, $\pi=0=\pi^{rev}$ a.e. in $\VV$, thus finally
$\pi=\pi^{rev}$ a.e. in $(0,T)\times\r^N$.\ep
\begin{Prop} Let $\pi\in\mbox{Vect}(\J^T)$. Then $\pi\in\R^T$
if and only if $|\pi|\in\W^T$.
\label{Prop abspi}\end{Prop}
\bp We first notice from \pref{piX}
that if $\pi\in\R^T$, then $|\pi|\in\R^T\subset\W^T$.
Conversely, if $|\pi|\in\W^T$, then by Proposition \ref{Lemma posVT},
$|\pi|\in\V^T$. By definition of $\V^T$, we deduce that $\pi\in\V^T$,
and by Proposition \ref{Prop jacVT}, $\pi\in\R^T$.\ep
	\begin{Rmk} \rm Any $\pi\in\W^T$ satisfies
$\pi\in C([0,T],L^1_{loc}(\r^N))$, and in particular
$|\pi|\in C([0,T],L^\infty_{loc}(\r^N)w*)$.
This is the consequence of Lemma \ref{Lemma OSLC} and of
Theorem 3.3 in \cite{B}.
	\label{Rmk strongcont}\end{Rmk}
	\begin{Prop} For any $\pi\in\W^T$, the application
\begin{equation}
	t \mapsto \| \pi (t, \cdot)\| _{L^1(\r^N)} \in [0,\infty]
	\label{tdonneJ}
\end{equation}
is nonincreasing in $[0,T]$.
Moreover,

\noindent (i) If $\pi$ is reversible then this function is constant.

\noindent (ii) If $\pi\in\mbox{Vect}(\J^T)$ and if the above
function is constant and finite, then $\pi$ is reversible.
	\label{Prop variation}\end{Prop}
\bp According to Lemma \ref{Lemma weakprod}, for any
$p \in \Lip_{loc}([0,T]\times\r^N)$ solving \pref{noncons},
we have $p\,\pi\in\W^T$. Take any transport flow $X^T$,
and choose $p(t,x)=\varphi(X^T(t,x))$,
with $\varphi\in C^1_c(\r^N)$.
Taking into account Remark \ref{Rmk proper} and
integrating the equation satisfied by $p\,\pi$,
we get for any $0\le t\le T$
\begin{equation}
	\int_{\r^N}\varphi(x)\pi(T,x)\,dx
	=\int_{\r^N}\varphi(X^T(t,x))\pi(t,x)\,dx.
	\label{ppiconst}
\end{equation}
This means that $\pi(T,.)dx$ is the image of $\pi(t,.)dx$ by
the function $x\mapsto X^T(t,x)$. In particular, \pref{ppiconst}
is valid also for all $\varphi$ measurable and
bounded with compact support. Then, taking $\varphi=\psi\sign(\pi(T,.))$,
we get for any $\psi$ measurable bounded
with compact support and nonnegative
\begin{equation}
	\int_{\r^N}\psi(x)|\pi(T,x)|\,dx
	\le\int_{\r^N}\psi(X^T(t,x))|\pi(t,x)|\,dx.
	\label{ppiabs}
\end{equation}
Letting $\psi\rightarrow 1$ we conclude that
\begin{equation}
	\int_{\r^N}|\pi(T,x)|\,dx
	\le\int_{\r^N}|\pi(t,x)|\,dx.
	\label{absdecr}
\end{equation}
Since the restriction on a subinterval of any weak solution $\pi$
to \pref{back} is again a weak solution, we can replace $T$
in \pref{absdecr} by any value $s\ge t$, which proves the first assertion.

In order to prove (i), we first recall
that if $\pi\in\R^T$, then $|\pi|\in\R^T$.
Thus in the case where $\pi(T,.)\in L^\infty_c(\r^N)$, we obviously
get that $\int|\pi(t,x)|dx$ is constant by integrating the equation
satisfied by $|\pi|$. In the general case $\pi(T,.)\in L^\infty_{loc}$,
we can define $\pi^T_n(x)=\pi(T,x)\1_{|x|\le n}$.
Denoting by $\pi_n$ the reversible solution with final data $\pi^T_n$,
we have that $\pi_n(t,.)\rightarrow \pi(t,.)$ a.e. for each $t$,
and we conclude by monotone convergence.

Let us finally prove (ii). Consider the reversible solution $\pi^{rev}$
with final data $\pi(T,.)$.
Since $\pi\in\mbox{Vect}(\J^T)$, we have by the same argument as
in the proof of Proposition \ref{Prop jacVT} that $\pi=\pi^{rev}$
a.e. in $(\VV)^c$. Thus for almost every $t$,
\begin{equation}
	\int_{\r^N}|\pi(t,x)|\,dx
	=\int_{(t,x)\notin\VV}|\pi^{rev}(t,x)|\,dx
	+\int_{(t,x)\in\VV}|\pi(t,x)|\,dx.	
	\label{decint}
\end{equation}
But since $\pi^{rev}=0$ a.e. in $\VV$, and by (i), we have
\begin{equation}
	\int_{(t,x)\notin\VV}|\pi^{rev}(t,x)|\,dx
	=\int_{\r^N}|\pi^{rev}(t,x)|\,dx
	=\int_{\r^N}|\pi(T,x)|\,dx.
	\label{eqpirev}
\end{equation}
Therefore, writing that $\int|\pi(t,x)|dx$ is constant and finite,
we obtain that $\int_{(t,x)\in\VV}|\pi(t,x)|\,dx=0$ for a.e. $t$, which
gives that $\pi=0$ a.e. in $\VV$. We conclude that
$\pi=\pi^{rev}$ a.e. in $(0,T)\times\r^N$.\ep

\Section{Forward problem, duality solutions}\label{forward}

We consider now the forward nonconservative transport problem
\pref{transport}-\pref{Cauchy}. We can only hope to have solutions
of bounded variation in $x$, thus we define the space
\begin{equation}
	\S=C([0,T],L^1_{loc}(\r^N))\cap\B([0,T],BV_{loc}(\r^N)),
	\label{SBV}
\end{equation}
where $\B$ stands for the space of bounded functions.
This regularity is however not enough to ensure a unique {\sl a priori}
determination of the product
$\aa\times\nabla u$. We therefore have to define solutions in a
weak sense, via the backward conservative problem studied in the
previous section. Reversible solutions on a subinterval $[0,\tau]$ are
involved, which is meaningful, $T$ has to be replaced by $\tau$
in Definitions \ref{rev-sol} and \ref{Def flow}.
	\begin{Def}[Duality solutions]\label{dual-sol}
We say that $u\in \S$ is a {\it duality solution} to \pref{transport} if
for any $0<\tau\le T$ and for any reversible solution $\pi$
to $\partial_t\pi+\diver(\aa\pi)=0$ in $(0,\tau)\times\r^N$
with compact support in $x$, one has that
\begin{equation}
	t\longmapsto \int_{\r^N} u(t,x)\, \pi (t,x)\, dx
	\ \mbox{ is constant in }[0,\tau].
	\label{defdual}
\end{equation}
	\end{Def}
This definition is motivated by the following lemma.
	\begin{Lm} Let $p\in\Lip_{loc}([0,T]\times\r^N)$ solve
$\partial_t p + \aa\cdot \nabla p =0$ a.e. in $(0,T)\times\r^N$.
Then $p$ is a duality solution.
	\label{Lemma lipdual}\end{Lm}
\bp By Lemma \ref{Lemma weakprod}, for any reversible solution
$\pi$ in $(0,T)\times\r^N$, we have
$\partial_t(p\pi)+\diver(\aa p\pi)=0$. If $\pi$ has compact support,
we deduce by integration in $x$ that $\int p(t,x)\pi(t,x)\,dx=cst$
in $[0,T]$. Since $p$ is also a solution on a subinterval $[0,\tau]$,
this yields the result.\ep

We are now in position to prove the main result of this section, namely existence
and uniqueness in the duality sense for the transport equation \pref{transport}.
	\begin{Th}[Nonconservative forward Cauchy problem]\label{exs-uniq}
For any $u^0\in BV_{loc}(\r^N)$, there exists a unique duality solution
$u \in \S$ to \pref{transport} such that $u(0,.)=u^0$.
	\end{Th}
{\bf Proof of uniqueness.}\;
Assume that there exists a duality solution $u$ such that $u(0,.)=0$.
Then, by definition, we have for any $0\le t\le \tau$
\begin{equation}
	\int_{\r^N} u(t,x) \, \pi(t,x) \, dx
	=\int_{\r^N} u(0,x) \, \pi(0,x) \, dx=0,
	\label{dualvanish}
\end{equation}
for any $\pi$ reversible solution in $[0,\tau]$ with compact support.
Choosing in particular $t=\tau$, and since $\pi(\tau,.)$ is arbitrary
in $L^\infty_c$, we obtain that $u(\tau,.)=0$.
This is true for any $0<\tau\le T$, thus $u\equiv 0$.\ep

\noindent The existence proof makes use of the following a priori $BV$ bound.
\begin{Lm}\label{Lemma BVbound}
Assume that $\aa\in C^1([0,T]\times\r^N)$ and that
$u\in C^2([0,T]\times\r^N)$ solves \pref{transport}
in the classical sense. Then for any $0\le t\le T$,
\begin{equation}
	\int\limits_{B(x_0,R)}\sum_{i=1}^N\left
	|\frac{\partial u}{\partial x_i}(t,x)\right|dx
	\le \sqrt{N}e^{(N-1)\int_0^t\alpha(s)ds}\mkern -20mu
	\int\limits_{B(x_0, R+t\|\aa\|_\infty)}\sum_{i=1}^N\left
	|\frac{\partial u}{\partial x_i}(0,x)\right|dx.
	\label{bvbound}
\end{equation}
\end{Lm}
\bp
Differentiating \pref{transport} with respect to $x_j$, we get
\begin{equation}
	(\partial_t+\aa\cdot\nabla)(\partial_j u)
	+\partial_j\aa\cdot\nabla u=0.
	\label{diu}
\end{equation}
Then, defining $\psi = ( \sum_j (\partial_j u )^2)^{1/2}$, we
obtain
\begin{equation}\label{equpsi}
	\psi\, (\partial_t + \aa \cdot \nabla) \psi
	+ \sum_{i,j} \partial_j \aa_i\,\partial_i u\, \partial_j u =0,
\end{equation}
or
\begin{equation}\label{equpsicons}
	\psi (\partial_t \psi + \diver(\aa \psi))
	+ \nabla u^t \nabla\aa\nabla u - \psi^2\diver\aa =0 .
\end{equation}
Noticing that $\nabla u^t\,\nabla u=\psi^2$ and
\begin{equation}
	\psi^2\diver\aa-\nabla u^t \nabla\aa\nabla u
	=\nabla u^t\left((\diver\aa)\I
	-\frac{\nabla\aa + \nabla\aa^t}{2}\right)\nabla u,
	\label{idennablaa}
\end{equation}
the equation \pref{equpsicons} gives with the inequality \pref{borna}
\begin{equation}
	\psi (\partial_t \psi + \diver (\aa \psi))
	\le (N-1)\,\alpha(t)\,\psi^2,
	\label{inegpsi}
\end{equation}
and therefore
\begin{equation}\label{vtpsi}
	\partial_t \psi + \diver (\aa \psi)
	\leq (N-1)\,\alpha(t)\,\psi.
\end{equation}
Defining $\phi=\psi e^{-(N-1)\int_0^t\alpha}$, we deduce that
\begin{equation}
	\partial_t \phi + \diver (\aa \phi)\leq 0,
	\label{inegphi}
\end{equation}
and integrating this over a cone yields
\begin{equation}
	\int\limits_{B(x_0,R)} \phi(t,x) \, dx
	\leq \int\limits_{B(x_0, R+t\|\aa\|_\infty)} \phi(0,x)\,dx.
	\label{estconephi}
\end{equation}
Translating this to $\psi$ gives
\begin{equation}
	\int\limits_{B(x_0,R)} \psi(t,x) \, dx
	\le e^{(N-1)\int_0^t\alpha}\mkern -20mu
	\int\limits_{B(x_0, R+t\|\aa\|_\infty)} \psi(0,x)\,dx.
	\label{estconepsi}
\end{equation}
But according to the Cauchy-Schwarz inequality
$\psi\le\sum|\partial_j u|\le \sqrt{N}\psi$,
which gives the estimate \pref{bvbound}.\ep
\begin{Rmk} \rm In the case $N=1$, Lemma \ref{Lemma BVbound}
reduces to a TVD property.
\label{Rmk TVD}\end{Rmk}
{\bf Proof of existence in Theorem \ref{exs-uniq}.}\;
Consider a sequence of coefficients $\aa_n\in C^2$, $\aa_n$
bounded in $L^\infty$, satisfying the OSLC condition for some $\alpha_n$
bounded in $L^1$, such that $\aa_n\to \aa$ in $L^1_{loc}$.
Then we have a classical flow $X_n(s,t,x)\in C^2$ associated to $\aa_n$,
satisfying $\partial_sX_n(s,t,x)=\aa_n(s,X_n(s,t,x))$, $X_n(t,t,x)=x$.
Consider also a sequence of initial data $u^0_n\in C^2$, bounded
in $BV_{loc}$, such that $u^0_n\rightarrow u^0$ in $L^1_{loc}$.
Then, let us define $u_n(t,x)\in C^2$ to be the classical solution to
\begin{equation}
	\partial_t u_n+\aa_n\cdot\nabla u_n=0\quad
	\mbox{in }(0,T)\times\r^N,\qquad u_n(0,.)=u_n^0.
	\label{transpun}
\end{equation}
According to Lemma \ref{Lemma BVbound},
$u_n$ is uniformly bounded in $\B([0,T],BV_{loc}(\r^N))$.
But since $u_n$ is given by $u_n(t,x)=u_n^0(X_n(0,t,x))$,
we have for any $0<\tau\le T$, $\varphi\in L^\infty_c(\r^N)$ and
$0\le t\le\tau$
\begin{equation}\begin{array}{l}
	\dsp\hphantom{=\ }\int_{\r^N}u_n(t,x)\varphi(X_n(\tau,t,x))J(X_n)(\tau,t,x)\,dx\\
	\dsp=\int_{\r^N}u_n(t,X_n(t,\tau,x))\varphi(x)\,dx\\
	\dsp=\int_{\r^N}u_n^0(X_n(0,\tau,x))\varphi(x)\,dx\\
	\dsp=\int_{\r^N}u_n^0(x)\varphi(X_n(\tau,0,x))J(X_n)(\tau,0,x)\,dx.
	\label{eqintun}
	\end{array}
\end{equation}
But as in the proof of Proposition \ref{Prop existflow}, $X_n$
is uniformly Lipschitz continuous in the domain $0\le t\le s\le T$,
$x\in\r^N$, and up to a subsequence it converges to $X$,
such that $X^\tau\equiv X(\tau,.,.)$ is a transport flow in $[0,\tau]$.
Taking $t=\tau$ in \pref{eqintun}, we get that
$|\int u_n(\tau,.)\varphi|\le C\|\varphi\|_\infty$, the constant
$C$ depending on the support of $\varphi$. Therefore
$u_n$ is bounded in $C([0,T],L^1_{loc}(\r^N))$ also.
But by \pref{transpun}, $u_n$ is equicontinuous in time, thus
extracting a subsequence again, $u_n\rightarrow u$ in
$C([0,T],L^1_{loc}(\r^N))$, with $u\in\S$ satisfying $u(0,.)=u^0$.
For any $\varphi\in C_c(\r^N)$ we can pass to the limit in \pref{eqintun},
thus for $0\le t\le\tau$
\begin{equation}\begin{array}{l}
	\dsp\hphantom{=\ } \int_{\r^N}u(t,x)\varphi(X(\tau,t,x))J(X)(\tau,t,x)\,dx\\
	\dsp =\int_{\r^N}u^0(x)\varphi(X(\tau,0,x))J(X)(\tau,0,x)\,dx.
	\label{eqintu}
	\end{array}
\end{equation}
By approximation this is still valid for $\varphi\in L^\infty_c$.
Noticing that $\pi(t,x)\equiv\varphi(X(\tau,t,x))J(X)(\tau,t,x)$ is the reversible
solution in $[0,\tau]$ with final data $\varphi$, we conclude that
\pref{defdual} holds, and $u$ is a duality solution.\ep

\Section{Weak stability}\label{stability}

To have weak stability results under sharp assumptions
is a key point in developing an efficient well-posedness theory.
We show here that our notions of reversible and duality
solutions are very well adapted to this achievement.\\

\noindent We shall consider in this section a sequence of coefficients
$\aa_n$ such that
\begin{equation}
	\aa_n \mbox{ is uniformly bounded in }L^\infty((0,T)\times\r^N),
	\label{borninf}
\end{equation}
and
\begin{equation}
	\aa_n \mbox{ satisfies an OSLC condition for some }\alpha_n
	\mbox{ bounded in }L^1(0,T).
	\label{bornoslc}
\end{equation}
Then, after extraction of a subsequence,
\begin{equation}
	\aa_n\rightharpoonup\aa\mbox{ in }L^\infty w*,
	\label{anlim}
\end{equation}
and according to Lemma \ref{Lemma suitean}, $\aa$ also satisfies
an OSLC condition.
The main two results of this section are the following.
	\begin{Th}[reversible backward stability] Assume \pref{borninf}-\pref{anlim}, and
let $\pi^T_n$ be a bounded sequence in $L^\infty_{loc}(\r^N)$
such that $\pi^T_n\rightharpoonup\pi^T$ locally in $L^\infty w*$.
Then the reversible solution $\pi_n$ to
\begin{equation}
	\partial_t\pi_n+\diver(\aa_n\pi_n)=0\mbox{ in }(0,T)\times\r^N,
	\qquad \pi_n(T,.)=\pi^T_n
	\label{pin}
\end{equation}
converges in $C([0,T],L^\infty_{loc}(\r^N)w*)$
to the reversible solution $\pi$ to
\begin{equation}
	\partial_t\pi+\diver(\aa\pi)=0\mbox{ in }(0,T)\times\r^N,
	\qquad \pi(T,.)=\pi^T.
	\label{pilim}
\end{equation}
Moreover, $\aa_n\pi_n\rightharpoonup\aa\pi$.
	\label{Th limrever}\end{Th}
	\begin{Th}[forward duality stability] Assume \pref{borninf}-\pref{anlim}, and
let $u^0_n$ be a bounded sequence in $BV_{loc}(\r^N)$
such that $u^0_n\rightarrow u^0$ in $L^1_{loc}$.
Then the duality solution $u_n$ to
\begin{equation}
	\partial_tu_n+\aa_n\cdot\nabla u_n=0\mbox{ in }(0,T)\times\r^N,
	\qquad u_n(0,.)=u^0_n
	\label{un}
\end{equation}
converges in $C([0,T],L^1_{loc}(\r^N))$ to the duality solution $u$ to
\begin{equation}
	\partial_tu+\aa\cdot\nabla u=0\mbox{ in }(0,T)\times\r^N,
	\qquad u(0,.)=u^0.
	\label{ulim}
\end{equation}
	\label{Th limdual}\end{Th}
{\bf Proof of Theorem \ref{Th limrever}.}\;
The sequence $\pi_n$ is bounded in $L^\infty_{loc}(\r^N)$,
uniformly in $t$ and $n$.
Since \pref{pin} gives compactness in time, we have after extraction
of a subsequence that $\pi_n$ converges to some $\pi$ in
$C([0,T],L^\infty_{loc}(\r^N)w*)$.
We have that
\begin{equation}
	\pi_n(t,x)=\pi_n^T(X_n^T(t,x))J(X_n^T)(t,x),
	\label{formulepin}
\end{equation}
where $X_n^T$ is a transport flow associated to $\aa_n$.
According to \pref{lipXn}, we can choose $X_n^T$ uniformly
bounded in $\Lip([0,T]\times\r^N)$, and thus we can extract
a subsequence converging locally uniformly
to some $X^T\in \Lip([0,T]\times\r^N)$.
Taking into account Lemma \ref{Lemma OSLC}, we have that
$\diver \aa_n$ is bounded in ${\cal M}_{loc}((0,T)\times\r^N)$.
Thus for $i=1,\dots,N$,
\begin{equation}\begin{array}{l}
	\dsp\hphantom{\longrightarrow} \aa_n\cdot\nabla (X_n^T)_{i}
	=\diver(\aa_n(X_n^T)_i)-\diver(\aa_n)(X_n^T)_i\\
	\dsp\longrightarrow \diver(\aa\,(X^T)_i)-\diver(\aa)(X^T)_i
	=\aa\cdot\nabla (X^T)_{i},
	\label{converXn}
	\end{array}
\end{equation}
and we deduce that $X^T$ is a transport flow associated to $\aa$.
According to Remark \ref{Rmk proper} and the uniform bounds,
when localizing in $x$ in \pref{formulepin}, only values of $\pi_n^T$
on a bounded set are involved. Therefore, we can apply
Lemma \ref{Lemma weakstab} to pass to the limit in \pref{formulepin}
at a fixed time $t$, and we conclude that
$\pi(t,x)=\pi^T(X^T(t,x))J(X^T)(t,x)$, i.e. that $\pi$
is the reversible solution to \pref{pilim}.
The uniqueness of the limit ensures that in fact it is not necessary
to extract any subsequence.

Let us finally prove the convergence of $\aa_n\pi_n$.
According to \pref{formulepin} and to \pref{Vkk}, we have
\begin{equation}
	(\aa_n)_k\pi_n=\pi_n^T(X_n^T)
	(-1)^k\det\left(\partial_tX_n^T,\frac{\partial X_n^T}{\partial x_1},
	\dots,\widehat{\frac{\partial X_n^T}{\partial x_{k}}},
	\dots,\frac{\partial X_n^T}{\partial x_{N}}\right).
	\label{anpin}
\end{equation}
Applying Lemma \ref{Lemma weakstab} in the variables $(t,x)$, we
conclude that we can pass to the limit weakly in the right-hand side,
which gives that $\aa_n\pi_n\rightharpoonup\aa\pi$.\ep

\begin{Rmk} \rm Another way of proving the convergence of $\aa_n\pi_n$ is by using
Lemma \ref{Lemma OSLC}. From the identity
$\Delta \aa_i=\sum_j\partial_j(\partial_j\aa_i+\partial_i\aa_j)
-\partial_i\diver\aa$ and the fact that $\nabla\aa_n+\nabla\aa_n^t$
is bounded in ${\cal M}_{loc}((0,T)\times\r^N)$, we deduce that
$\aa_n$ is compact in $x$, i.e. that
$\|\aa_n(t,x+h)-\aa_n(t,x)\|_{L^1((0,T)\times B_R)}\rightarrow 0$
as $h$ tends to $0$, uniformly in $n$.
Since $\pi_n\rightarrow\pi$ in $C([0,T],L^\infty_{loc}(\r^N)w*)$,
this is enough to conclude that $\aa_n\pi_n\rightharpoonup\aa\pi$.
\end{Rmk}

\begin{Rmk} \rm According to Lemma \ref{Lemma OSLC} and to
Theorem 3.3 in \cite{B}, $\pi_n$ and $\pi$ lie in $C([0,T],L^1_{loc}(\r^N))$.
It is an open problem to prove that if $\pi_n^T\rightarrow\pi^T$
in $L^1_{loc}$, then $\pi_n\rightarrow \pi$ in $C([0,T],L^1_{loc}(\r^N))$.
\label{Rmk strongCV}\end{Rmk}
\noindent{\bf Proof of Theorem \ref{Th limdual}.}\;
The sequence $u_n$ is bounded in $\S$, and equicontinuous in time
(even if \pref{un} does not hold in the classical sense, the
a priori estimate on $\partial_t u_n$ in $L^\infty_t({\cal M}_{loc}(\r^N))$
is valid). Therefore, after extraction of a subsequence,
$u_n$ converges in $C([0,T],L^1_{loc}(\r^N))$ to some $u\in\S$
satisfying $u(0,.)=u^0$. Now, for any $0<\tau\le T$
and any $\pi^\tau\in L^\infty_c(\r^N)$, we have that
\begin{equation}
	t\longmapsto \int_{\r^N} u_n(t,x)\, \pi_n (t,x)\, dx
	\ \mbox{ is constant in }[0,\tau],
	\label{defdualn}
\end{equation}
where $\pi_n$ is the reversible solution to
$\partial_t\pi_n+\diver(\aa_n\pi_n)=0$ in $(0,\tau)\times\r^N$,
$\pi_n(\tau,.)=\pi^\tau$. Applying Theorem \ref{Th limrever}
on $(0,\tau)$, we deduce that $\pi_n\rightarrow \pi$
in $C([0,\tau],L^\infty_{loc}(\r^N)w*)$ where $\pi$ is the
reversible solution to
$\partial_t\pi+\diver(\aa\pi)=0$ in $(0,\tau)\times\r^N$,
$\pi(\tau,.)=\pi^\tau$.
Passing to the limit in \pref{defdualn}, we get that
\begin{equation}
	t\longmapsto \int_{\r^N} u(t,x)\, \pi (t,x)\, dx
	\ \mbox{ is constant in }[0,\tau],
	\label{defduallim}
\end{equation}
which means that $u$ is the duality solution to \pref{ulim}.
By uniqueness, we conclude that it is indeed not necessary to extract
any subsequence.\ep

\Section{Nonuniqueness of transport flows}\label{nonuniqueness}

The aim of this section is to show that the uniqueness results
established in Section \ref{BackProblem} are optimal, in
the sense that first the transport flow is not unique,
and second its jacobian is uniquely determined.\\

Our example is in two space dimensions, and is really the simplest
nontrivial coefficient we can think of,
\begin{equation}
	\aa(t,x_1,x_2)=(-\sign x_1,0).
	\label{asigne}
\end{equation}
It is bounded, and obviously satisfies the OSLC condition
with $\alpha\equiv 0$. Since the second component of $\aa$ vanishes
identically, the second variable $x_2$ only
stands as a parameter in \pref{back} and \pref{noncons}.
We deduce the general solution $p\in \Lip_{loc}([0,T]\times\r^N)$
to \pref{noncons} by just adding a parameter $x_2$ to the
general solution of the one-dimensional problem obtained in \cite{BJ},
\begin{equation}
	p(t,x_1,x_2)=\varphi\Bigl(\bigl(|x_1|-(T-t)\bigr)_+\sign x_1,x_2\Bigr)
	+h\Bigl(\bigl(T-t-|x_1|\bigr)_+,x_2\Bigr),
	\label{solp}
\end{equation}
with
\begin{equation}
	\varphi\in \Lip_{loc}(\r\times\r),\qquad
	h\in\Lip_{loc}([0,T]\times\r),\ h(0,x_2)=0.
	\label{phih}
\end{equation}
Similarly, the general solution $\pi\in C([0,T],L^\infty_{loc}(\r^N)w*)$
to \pref{back} (i.e. $\pi\in\W^T$), is given by
\begin{equation}\begin{array}{l}
	\dsp \pi(t,x_1,x_2)=\1_{|x_1|\ge T-t}
	\,\psi\Bigl(\bigl(|x_1|-(T-t)\bigr)_+\sign x_1,x_2\Bigr)\\
	\dsp\hphantom{\pi(t,x_1,x_2)=}
	+\1_{|x_1|<T-t}\,g\Bigl(\bigl(T-t-|x_1|\bigr)_+,x_2\Bigr) \sign x_1,
	\label{solpi}
	\end{array}
\end{equation}
with
\begin{equation}
	\psi\in L^\infty_{loc}(\r\times\r),\qquad
	g\in L^\infty_{loc}([0,T]\times\r).
	\label{psig}
\end{equation}
Now, we see from \pref{solp} that $p\in\E^T$ if and only if $\varphi\equiv 0$,
and therefore
\begin{equation}
	\VV=\{ (t,x_1,x_2)\in (0,T)\times\r^2
	\ \mbox{such that}\ |x_1|< T-t \}.
	\label{VVT}
\end{equation}
Then, we see that $\pi$ in \pref{solpi} lies in $\V^T$ if and only if
$g\equiv 0$.\\

Let us now look for transport flows $X^T$, as in Definition \ref{Def flow}.
Each component of $X^T$ has to be of the form \pref{solp},
with final data $\varphi_i$ given respectively by $\varphi_1(x_1,x_2)=x_1$
and $\varphi_2(x_1,x_2)=x_2$. The functions $h_i$ are arbitrary, thus
\begin{equation}\begin{array}{l}
	\dsp X^T(t,x_1,x_2)=\Biggl(\bigl(|x_1|-(T-t)\bigr)_+\sign x_1
	+h_1\Bigl(\bigl(T-t-|x_1|\bigr)_+,x_2\Bigr),\Biggr.\\
	\dsp\hphantom{X^T(t,x_2,x_2)=\qquad} \Biggl.x_2
	+h_2\Bigl(\bigl(T-t-|x_1|\bigr)_+,x_2\Bigr)\Biggr),
	\label{XXT}
	\end{array}
\end{equation}
and we need only to write that $J(X^T)\ge 0$. According to the
conditions \pref{phih}, in the set where
$|x_1|\ge T-t$ we have $J(X^T)=1$, while
$J(X^T)=-(\sign x_1)\left(\partial_t h_1(1+\partial_2h_2)-\partial_2h_1\partial_th_2\right)(\bigl(T-t-|x_1|\bigr)_+,x_2)$
in the set $|x_1|< T-t$. Since this last quantity is odd with respect
to $x_1$, its nonnegativity implies that it vanishes identically, thus
the condition $J(X^T)\ge 0$ resumes to
\begin{equation}
	\partial_t h_1\,(1+\partial_2h_2)-\partial_2h_1\,\partial_th_2\equiv 0.
	\label{jac0}
\end{equation}
Once this is satisfied, we have
\begin{equation}
	J(X^T)=\1_{|x_1|\ge T-t},
	\label{jacXT}
\end{equation}
and as predicted by Lemma \ref{Lemma uniqJ}, this is independent
of the choice of $h_1$, $h_2$. According to Definition \ref{rev-sol},
$\pi$ in \pref{solpi} is reversible if and only if $g\equiv 0$.
Thus for our coefficient $\aa$, we have here equality in \pref{inclusion}.\\

Finally, we can observe that the transport flow is not unique even
if we impose the semi-group property. Denoting by $X(s,t,x)=X^s(t,x)$
for $0\le t\le s\le T$, this means that
\begin{equation}
	X(s,t,X(t,\tau,x))=X(s,\tau,x)\quad
	\mbox{for}\quad 0\le \tau\le t\le s\le T.
	\label{semig}
\end{equation}
Choosing $h_1\equiv 0$ and $h_2(t,x_2)=\lambda t$ in \pref{XXT},
the conditions \pref{phih} and \pref{jac0} are satisfied thus
we get a transport flow $X_\lambda$ for any $\lambda\in\r$,
\begin{equation}
	X_\lambda(s,t,x_1,x_2)=\Biggl(\bigl(|x_1|-(s-t)\bigr)_+\sign x_1,
	x_2+\lambda\bigl(s-t-|x_1|\bigr)_+\Biggr).
	\label{XXlambda}
\end{equation}
Noticing that $(\tau_1-|x_1|)_++(\tau_2-(|x_1|-\tau_1)_+)_+=(\tau_1+\tau_2-|x_1|)_+$
for any $\tau_1,\tau_2\ge 0$, one easily checks that \pref{semig}
holds for any $\lambda$. However, it is well-known that for a coefficient
satisfying the OSLC condition, there is only one flow $X$ which
solves $\partial_s X(s,t,x)=\aa(s,X(s,t,x))$ in the sense of Filippov.
It is indeed obtained for the choice $\lambda=0$.


\def\theequation{A.\arabic{equation}}
\def\thesection{A}
\def\thesubsection{A}
\def\thesubsubsection{A}
\setcounter{Th}{0}
\setcounter{Rmk}{0}
\setcounter{Ex}{0}
\setcounter{equation}{0}
\section*{Appendix}\label{appendix}
This appendix is devoted to some useful results on jacobian determinants.
The first two lemmas are proved in \cite{B}, respectively Lemma 2.5 and Theorem 2.4.
	\begin{Lm}[Divergence chain rule] Let $g\in C^1(\r^d,\r^d)$ such that
\begin{equation}
	\diver g\in L^\infty,\qquad |g(y)|\le C(1+|y|).
	\label{diverg}
\end{equation}
Let $\Omega$ be an open subset of $\r^N$, and
$u\in L^1_{loc}(\Omega,\r^d)$, $d\le N$, such that
$\partial_{x_1}u,\cdots,\partial_{x_N}u\in L^d_{loc}(\Omega,\r^d)$.
Then for any injective $\sigma:\N_d\rightarrow\N_N$,
\begin{equation}
	\sum_{k=1}^d\frac{\partial}{\partial x_{\sigma(k)}}
	\left[I\left(\left[\frac{\partial u_i}{\partial x_{\sigma(j)}}
	\right]_{1\le i,j\le d}\right)
	g(u)\right]_k=(\diver g)\circ u\ \det\left[
	\frac{\partial u_i}{\partial x_{\sigma(j)}}
	\right]_{1\le i,j\le d},
	\label{divchainrule}
\end{equation}
where $I(A)$ denotes the pseudo inverse (i.e. transpose of the cofactor matrix)
of $A\in{\mathbb M}_d(\r)$, $I(A)=(\mbox{com}A)^t$.
	\label{Lemma formulediv}\end{Lm}
	\begin{Lm} Let $\Omega$ be an open subset of $\r^N$ and
$u\in L^1_{loc}(\Omega,\r^d)$, $d\le N$, such that
$\partial_{x_1}u,\cdots,\partial_{x_N}u\in L^d_{loc}(\Omega,\r^d)$.
Then for any Borel set $Z\subset\r^d$ such that $|Z|=0$ and any injective
$\sigma:\N_d\rightarrow\N_N$,
\begin{equation}
	\left|\left\{x\in\Omega\ ;\ u(x)\in Z\mbox{ and }
	\det\left(\frac{\partial u_i}{\partial x_{\sigma(j)}}
	\right)_{1\le i,j\le d}\not=0\right\}\right|=0.
	\label{measurezero}
\end{equation}
In other words, $\mathop{\rm{rank}}(Du(x))<d$ almost everywhere
on $u^{-1}(Z)$.
	\label{Lemma detnulae}\end{Lm}
	\begin{Lm} Let $\Omega$ be an open subset of $\r^{N+1}$ with $N\ge 1$, and
$H_1,\dots,H_N\in L^1_{loc}(\Omega)$ such that $\nabla H_i\in L^N_{loc}(\Omega)$.
Define for $k=0,\dots,N$
\begin{equation}
	V_k=(-1)^k\det\left(\frac{\partial H_i}{\partial x_j}\right)
	_{i=1,\dots,N,\ j=0,\dots,k-1,k+1,\dots,N}.
	\label{Ik}
\end{equation}
Then $V\in L^1_{loc}(\Omega)$ satisfies $\diver V=0$ in the sense
of distributions in $\Omega$.
	\label{Lemma divnul}\end{Lm}
\bp Let us assume that $H_i\in C^\infty$, the general case being
easily deduced by approximation. Consider an arbitrary function
$H_0\in C^\infty$, and apply Lemma \ref{Lemma formulediv} with $d=N+1$,
$g(y)=(1,0,\dots,0)\in\r^{N+1}$, $u=(H_0,\dots,H_N)$, and $\sigma=\I$.
We get that $\diver(I(\nabla u)g)=0$, which gives the result
since $V=(-1)^N I(\nabla u)g$.\ep

\noindent We can notice in the previous lemma that $H_i$ are like
Hamiltonians for $V$ since $V\cdot\nabla H_i=0$.
This can be seen by the identity
\begin{equation}\begin{array}{l}
	\dsp V\cdot\nabla H_i
	=(-1)^N\left[I(\nabla u)g\right]\cdot\nabla H_i\\
	\dsp\hphantom{V\cdot\nabla H_i}
	=(-1)^N\det\left(\nabla H_1,
	\dots,\nabla H_N,\nabla H_i\right)=0.
	\label{hamili}
	\end{array}
\end{equation}
\noindent Another proof of Lemma \ref{Lemma divnul} is to consider
the $N$ differential form $\omega=dH_1\wedge \dots\wedge dH_N$.
Then its external differential vanishes, $d\omega=0$, because $d^2=0$.
One can check that
$\omega=\sum_{k=0}^N(-1)^{k-1}V_k\ dx_0\wedge\dots\wedge
\widehat{dx_k}\wedge\dots\wedge dx_N$, giving $\diver V=0$.\\

\noindent The next lemma is a generalization of the weak stability
of the jacobian determinant (see \cite{Resh,Resh2} and \cite{Ball}).
	\begin{Lm} Let $\Omega$ be an open subset of $\r^N$ and $1\le d\le N$.
Consider a sequence $u^{(n)}$ bounded in $L^1_{loc}(\Omega,\r^d)$
such that $\partial_{x_j}u^{(n)}$ is bounded in $L^d_{loc}(\Omega,\r^d)$
for $1\le j\le N$, and assume that $u^{(n)}\rightarrow u$
in $L^1_{loc}(\Omega,\r^d)$, with $\partial_{x_j}u\in L^d_{loc}(\Omega,\r^d)$.
Consider also a sequence $\psi_n$ bounded in $L^\infty(\r^d)$
such that $\psi_n\rightharpoonup \psi$ in $L^\infty w*$. Then
for any $\sigma:\N_d\rightarrow \N_N$ injective,
\begin{equation}
	\psi_n(u^{(n)})\det\left[\frac{\partial u_i^{(n)}}
	{\partial x_{\sigma(j)}}\right]_{1\le i,j\le d}
	\rightharpoonup\qquad
	\psi(u)\det\left[\frac{\partial u_i}
	{\partial x_{\sigma(j)}}\right]_{1\le i,j\le d}
	\label{limdet}
\end{equation}
in the sense of distributions in $\Omega$.	
	\label{Lemma weakstab}\end{Lm}
\bp According to Lemma \ref{Lemma detnulae}, both sides
in \pref{limdet} are well-defined as elements of $L^1_{loc}(\Omega)$.

Let us consider first the case $\psi_n=\psi=1$, i.e.
the case of pure jacobians. Then the result can be established
by induction on $d$. When $d=1$, this is obvious by linearity.
Assuming the result to be true at level $d-1$, with $2\le d\le N$, we
apply Lemma \ref{Lemma formulediv} with $g(y)=y/d$, to $u^{(n)}$
and $u$. We deduce that we only need to prove that
\begin{equation}
	I\left(A_n\right)u^{(n)}
	\quad\rightharpoonup\quad
	I\left(A\right)u,
	\label{limprodI}
\end{equation}
where
\begin{equation}
	A_n=\left[\frac{\partial u_i^{(n)}}{\partial x_{\sigma(j)}}
	\right]_{1\le i,j\le d},\qquad
	A=\left[\frac{\partial u_i}{\partial x_{\sigma(j)}}
	\right]_{1\le i,j\le d}.
	\label{An}
\end{equation}
Since the coefficients of $I(A_n)$ are sub-determinants of
order $d-1$ of $A_n$, the recurrence assumption ensures that
$I(A_n)\rightharpoonup I(A)$. Since $I(A_n)$ is bounded
in $L^{d'}_{loc}$, this convergence holds in $L^{d'}_{loc}-w*$.
But according to Sobolev imbedding, $u^{(n)}\rightarrow u$
in $L^d_{loc}$, and this yields \pref{limprodI}.

In the case $\psi\not \equiv 1$, define
\begin{equation}
	g(y)=\frac{1}{|S^{d-1}|}\int_{\r^d}\left(\frac{y-z}{|y-z|^d}
	+\frac{z}{|z|^d}
	-\1_{|z|\ge 1}\left(\frac{y}{|z|^d}
	-d\times (y\cdot z)\frac{z}{|z|^{d+2}}\right)\right)\psi(z)\,dz.
	\label{convol}
\end{equation}
Then $g\in C^{0,\beta}_{loc}(\r^d,\r^d)$ for any $\beta<1$,
$\diver g=\psi$, and
\begin{equation}
	|g(y)|\le C_d(1+|y|\ln^+|y|)\|\psi\|_\infty.
	\label{croissg}
\end{equation}
This last estimate can be obtained as follows. For $R>0$, define
$M_R$ to be the supremum over $|y|\le R$ of the $L^1$ norm in $z$
of the term between parentheses in \pref{convol}. Then one can check
easily that $M_{2R}\le 2M_R+C_d R$, and this implies that
$M_R\le C_d(1+R\ln^+R)$. This proves \pref{croissg}.
Now, we claim that \pref{divchainrule} holds with this nonlinearity $g$.
Indeed the right-hand side makes sense according to
Lemma \ref{Lemma detnulae}, and the left-hand side also by \pref{croissg},
since by Sobolev imbedding, $u\in L^{d+\varepsilon}_{loc}$ for
some $\varepsilon>0$. The validity of \pref{divchainrule} for
$g$ defined by \pref{convol} and $\psi\in L^\infty(\r^d)$
can be obtained by approximating $\psi$
by some smooth functions with compact support.

\noindent Finally, for our stability result when $\psi_n\rightharpoonup\psi$
in $L^\infty w*$, according to the previous argument, we
can apply \pref{divchainrule} to $u$ and $g$ corresponding to $\psi$,
and also to $u^{(n)}$ and $g_n$ corresponding to $\psi_n$.
Thus in order to get \pref{limdet}, it is enough to prove that
$I(A_n)g_n(u^{(n)})\rightharpoonup I(A)g(u)$, with $A_n$ and $A$
defined by \pref{An}. By the first step of pure jacobians, we
have $I(A_n)\rightharpoonup I(A)$ in $L^{d'}_{loc}-w*$.
But since $\psi_n\rightharpoonup\psi$ we have $g_n\rightarrow g$
locally uniformly. Using \pref{croissg} and the fact that
$u^{(n)}\rightarrow u$ in $L^{d+\varepsilon}_{loc}$
for some $\varepsilon>0$, we get that
$g_n(u^{(n)})\rightarrow g(u)$ in $L^{d}_{loc}$,
and thus $I(A_n)g_n(u^{(n)})\rightharpoonup I(A)g(u)$.\ep

\end{document}